\documentclass[10pt]{article}

\usepackage{amsmath}
\usepackage{amssymb}
\usepackage{indentfirst}
\usepackage{graphics} 
\usepackage{color}

\makeatletter

\@addtoreset{equation}{section}
\makeatother
\def\<{\langle}
\def\>{\rangle}

\newtheorem{lem}{Lemma}[section]
\newtheorem{theo}{Theorem}[section]
\newtheorem{rem}{Remark}[section]
\newtheorem{pro}{Proposition}[section]

\makeatletter
   
   \@addtoreset{equation}{section}
\makeatother

\begin{document}
\title{\bf $L^{2}$-blowup estimates of the plate equation}

\author{Ryo Ikehata\thanks{ikehatar@hiroshima-u.ac.jp} \\ {\small Department of Mathematics}\\ {\small Division of Educational Sciences}\\ {\small Graduate School of Humanities and Social Sciences} \\ {\small Hiroshima University} \\ {\small Higashi-Hiroshima 739-8524, Japan}}
\date{}
\maketitle
\begin{abstract}
We consider the Cauchy problems in ${\bf R}^{n}$ for the plate equation with a weighted $L^{1}$-initial data. We derive optimal estimates of the $L^{2}$-norm of solutions for $n = 1,2,3,4$. In particular, such obtained results express infinite time blowup property in the case when the $0$-th moment of the initial velocity does not vanish. The idea to derive them is strongly inspired from an already developed technique \cite{I-14, IO, Ikehata}.  
\end{abstract}
\section{Introduction}
\footnote[0]{Keywords and Phrases: Plate equation; weighted $L^{1}$-data; low dimension; blowup in infinite time; sharp estimates.}
\footnote[0]{2010 Mathematics Subject Classification. Primary 35L05; Secondary 35B40, 35C20.}
This is one more application to the plate and/or beam equation of the method recently applied in \cite{Ikehata} to the wave equation. A similar arugument to these wave-like estimates to capture a singularity is introduced in \cite{I-14, IO}, and their related papers are published by the author's collaborative works (see \cite{I-10} and the references therein).\\ 
Now, we consider the Cauchy problem of the plate equation:
\begin{align}
& u_{tt} +\Delta^{2}u = 0,\ \ \ (t,x)\in (0,\infty)\times {\bf R}^{n},\label{eqn}\\
& u(0,x)= u_0(x), \quad  u_{t}(0,x)= u_{1}(x),\ x\in{\bf R}^{n}.\label{initial}
\end{align}
Here, we assume, for the moment, $[u_{0},u_{1}] \in H^{2}({\bf R}^{n})\times L^{2}({\bf R}^{n})$.\\ 
\noindent
Concerning the existence of a unique energy solution to problem \eqref{eqn}-\eqref{initial}, by the standard semi-group theory one can find that the problem (1.1)-(1.2) has a unique weak solution
\[u \in C([0,\infty);H^{2}({\bf R}^{n})) \cap C^{1}([0,\infty);L^{2}({\bf R}^{n})) \]
satisfying the energy conservation law such that
\begin{equation}\label{i-1}
E(t) = E(0),\quad t \geq 0,
\end{equation}
where the total energy $E(t)$ for the solution to problem \eqref{eqn}-\eqref{initial} can be defined by
\[E(t) := \frac{1}{2}\left(\Vert u_{t}(t,\cdot)\Vert^{2} + \Vert\Delta u(t,\cdot)\Vert^{2}\right).\] 
Here, $\Vert u\Vert$ denotes the usual $L^{2}$-norm of $u \in L^{2}({\bf R}^{n})$.

As is frequently appeared in several research papers on the plate and/or beam models, the assumption $n \geq 5$ is sometimes imposed to treat the problems in unbounded domains, and for this observation one can cite several typical papers such as \cite{D}, \cite{LC}, \cite{SJ}, \cite{L}, \cite{lin}, \cite{Miao} and \cite{St-2}, however, it seems that there are no any related papers to investigate actively a reason why $n \geq 5$ (see also \cite{peral} for a topic on the $L^{p}$-regularity condition). The purpose of this paper is to give a partial answer on this kind of problem. In this connection, in \cite{P, Racke} several decay estimates for the quantity $\Vert(u_{t}(t,\cdot),\Delta u(t,\cdot))\Vert_{q}$ can be studied in order to apply them to nonlinear problems. Our main concern is in observing a singularity near small frequency region, so we are particular about dealing with the asymptotic behavior of the quantity $\Vert u(t,\cdot)\Vert$. By the way, it should be mentioned that the equation considered in \cite{L} does not have any singularity near $0$ frequency region, so a topic taken up in \cite{L} seems to be a little different from ours.\\  
 
Before going to introduce our Theorems, we set
\[I_{0,n} := \Vert u_{0}\Vert + \Vert u_{0}\Vert_{L^{1}({\bf R}^{n})} +\Vert u_{1}\Vert + \Vert u_{1}\Vert_{L^{1}({\bf R}^{n})}.\]
The following two theorems give a hint about a question why $n \geq 5$ in the plate equation.
\begin{theo}\label{theorem1}
Let $n = 1,2,3$. Let $[u_{0},u_{1}] \in H^{2}({\bf R}^{n}) \times L^{2}({\bf R}^{n})$. Then, the solution $u(t,x)$ to problem \eqref{eqn}-\eqref{initial} satisfies the following properties under the additional regularity on the initial data:
\[[u_{0},u_{1}] \in L^{1}({\bf R}^{n}) \times L^{1}({\bf R}^{n}) \quad \Rightarrow \quad \Vert u(t,\cdot)\Vert \leq C_{1}I_{0,n}t^{1-\frac{n}{4}},\]
\[[u_{0},u_{1}] \in L^{1}({\bf R}^{n}) \times L^{1,1}({\bf R}^{n}) \quad \Rightarrow \quad C_{2}\left\vert\int_{{\bf R}^{n}}u_{1}(x)dx\right\vert t^{1-\frac{n}{4}} \leq \Vert u(t,\cdot)\Vert\]
for $t \gg 1$, where $C_{j} > 0$ {\rm ($j = 1,2$)} are constants depending only on $n$.
\end{theo}
Our next result is the case of $n = 4$.
\begin{theo}\label{theorem2}
Let $n = 4$. Let $[u_{0},u_{1}] \in H^{2}({\bf R}^{4}) \times L^{2}({\bf R}^{4})$. Then, the solution $u(t,x)$ to problem \eqref{eqn}-\eqref{initial} satisfies the following properties under the additional regularity on the initial data:
\[[u_{0},u_{1}] \in L^{1}({\bf R}^{4}) \times L^{1}({\bf R}^{4}) \quad \Rightarrow \quad \Vert u(t,\cdot)\Vert \leq C_{1}I_{0,4}\sqrt{\log t},\]
\[[u_{0},u_{1}] \in L^{1}({\bf R}^{4}) \times L^{1,1}({\bf R}^{4}) \quad \Rightarrow \quad C_{2}\left\vert\int_{{\bf R}^{4}}u_{1}(x)dx\right\vert\sqrt{\log t} \leq \Vert u(t,\cdot)\Vert\]
for $t \gg 1$, where $C_{j} > 0$ {\rm ($j = 1,2$)} are constants depending only on the space dimension $n = 4$. 
\end{theo}
\begin{rem}{\rm One of our advantages in the results is that we never use the compactness argument as is frequently developed in the wave equation case because of the non-Kowalewskian property of the equation \eqref{eqn}. The method will be widely applicable to the other types of evolution equations. Indeed, it will be possible to generalize our results to the more general $\sigma$-evolution equations:
\[u_{tt} + (-\Delta)^{\sigma}u = 0,\]
where $((-\Delta)^{\sigma}f)(x) := {\cal F}_{\xi \to x}^{-1}(\vert\xi\vert^{2\sigma}\hat{f}(\xi))(x)$ and $\sigma > 0$. A threshold number $n^{*}$ on the dimension $n$ to divide whether the blow up phenomenon occurs or not can be defined by $n^{*} = 2\sigma$.}
\end{rem}
\begin{rem}{\rm In both Theorems, the so-called infinite time blowup property can be observed in the case when $n < 5$. One has to treat more delicately when we study decay estimates, asymptotic profiles, nonlinear problems and so on, of the plate equation considered in unbounded domains, and in particular, in ${\bf R}^{n}$.}
\end{rem}

Contrary to infinite time blowup results above, in the case when $\int_{{\bf R}^{n}}u_{1}(x)dx = 0$, one has the $L^{2}$-boundedness property (at least) for $n = 3,4$.

\begin{theo}\label{theorem3}
Let $n = 3, 4$, and $[u_{0},u_{1}] \in H^{2}({\bf R}^{n}) \times L^{2}({\bf R}^{n})$. Assume further that $u_{1} \in L^{1,\gamma}({\bf R}^{n})$ with  
\[n = 4 \quad \Rightarrow \quad \gamma \in (0,1],\]
\[n = 3 \quad \Rightarrow \quad \gamma \in (\frac{1}{2},1].\]
If
\[\int_{{\bf R}^{n}}u_{1}(x)dx = 0,\]
then the solution $u(t,x)$ to problem \eqref{eqn}-\eqref{initial} satisfies 
\[\Vert u(t,\cdot)\Vert \leq C(\Vert u_{0}\Vert + \Vert u_{1}\Vert +  \Vert u_{1}\Vert_{1,\gamma}),\]
where $C > 0$ is a constant depending on the space dimension $n$ and $\gamma$. 
\end{theo}
In the case of $n = 1,2$, if one assumes an additional vanishing condition of the $1$th-order moment of the initial velocity $u_{1}$ one can state the following $L^{2}$-boundedness result.
\begin{theo}\label{theorem4}
Let $n = 1, 2$, and $[u_{0},u_{1}] \in H^{2}({\bf R}^{n}) \times L^{2}({\bf R}^{n})$. Assume further that $u_{1} \in L^{1,\gamma}({\bf R}^{n})$ with  
\[n = 2 \quad \Rightarrow \quad \gamma \in (1,2],\]
\[n = 1 \quad \Rightarrow \quad \gamma \in (\frac{3}{2},2].\]
Under the following two conditions,  
\begin{equation}\label{1000}
\int_{{\bf R}^{n}}u_{1}(x)dx = 0,\qquad \int_{{\bf R}^{n}}x_{j}u_{1}(x)dx = 0\quad (j = 1,2,\cdots, n),
\end{equation}
the solution $u(t,x)$ to problem \eqref{eqn}-\eqref{initial} satisfies 
\[\Vert u(t,\cdot)\Vert \leq C(\Vert u_{0}\Vert + \Vert u_{1}\Vert +  \Vert u_{1}\Vert_{1,\gamma}),\]
where $C > 0$ is a constant depending on the space dimension $n$ and $\gamma$. 
\end{theo}
{\bf Example.}\,Let $n = 1$, and let $g \in C_{0}^{\infty}({\bf R})$ be an odd function, and choose $u_{1}(x) := g'(x) \in C_{0}^{\infty}({\bf R})$. Then, one has
\[\int_{{\bf R}}u_{1}(x)dx = 0,\quad \int_{{\bf R}}x u_{1}(x)dx = \int_{{\bf R}}x g'(x)dx = -\int_{{\bf R}}g(x)dx = 0.\]
This expresses an example of the initial velocity satisfying \eqref{1000} for $n = 1$.\\
Similarly, in the case when $n = 2$, one can construct an example satsfying \eqref{1000} by
\[u_{1}(x_{1},x_{2}) := g'(x_{1})h'(x_{2}),\]
where $g,h \in C_{0}^{\infty}({\bf R})$.\\
As a counter part of Theorem 1.4 one can get the infinite time blowup result in the case when $n = 1,2$. For simplicity, we mention only the case of $\gamma = 2$ in order to compare it with Theorem 1.4..
\begin{theo}\label{theorem5}
{\bf (1)}\,Let $n = 1$, and $[u_{0},u_{1}] \in H^{2}({\bf R}) \times L^{2}({\bf R})$. Assume further that $u_{1} \in L^{1,2}({\bf R})$ satisfies  
\begin{equation}\label{1001}
\left\vert \int_{{\bf R}}u_{1}(x)dx\right\vert + \left\vert \int_{{\bf R}}xu_{1}(x)dx\right\vert > 0.
\end{equation}
Then the solution $u(t,x)$ to problem \eqref{eqn}-\eqref{initial} satisfies 
\[C\left(\left\vert\int_{{\bf R}}x u_{1}(x)dx\right\vert t^{\frac{1}{4}} + \left\vert\int_{{\bf R}}u_{1}(x)dx\right\vert t^{\frac{3}{4}}\right) \leq \Vert u(t,\cdot)\Vert,\quad t \gg 1,\]
where $C > 0$ is a constant. \\
{\bf (2)}\,Let $n = 2$, and $[u_{0},u_{1}] \in H^{2}({\bf R}^{2}) \times L^{2}({\bf R}^{2})$. Assume further that $u_{1} \in L^{1,2}({\bf R}^{2})$ satisfies  
\begin{equation}\label{1001}
\left\vert \int_{{\bf R}^{2}}u_{1}(x)dx\right\vert + \left\vert \int_{{\bf R}^{2}}x_{1}u_{1}(x)dx\right\vert + \left\vert \int_{{\bf R}^{2}}x_{2}u_{1}(x)dx\right\vert   > 0.
\end{equation}
Then the solution $u(t,x)$ to problem \eqref{eqn}-\eqref{initial} satisfies 
\[C\left(\left\vert\int_{{\bf R}^{2}}x_{1}u_{1}(x)dx\right\vert + \left\vert\int_{{\bf R}^{2}}x_{2}u_{1}(x)dx\right\vert\right)\sqrt{\log t} + C\left\vert \int_{{\bf R}^{2}}u_{1}(x)dx\right\vert\sqrt{t} \leq \Vert u(t,\cdot)\Vert,\quad t \gg 1,\]
where $C > 0$ is a constant. \\
\end{theo}
\begin{rem}{\rm As a result of Theorem 1.5, (for examle) in case of $n=1$,  even if $\int_{{\bf R}}u_{1}(x)dx = 0$, if $\int_{{\bf R}}x u_{1}(x)dx \ne 0$, then one has the infinite time blowup result with its rate $t^{1/4}$, however one can know nothing about the optimality of the blowup rate $t^{1/4}$. A corresponding result for $n = 2$ is true similarly. Of course, one should investigate (optimal) upper bound estimates of Theorem 1.5, however, this study will be left to the reader's interest.}
\end{rem}

{\bf Notation.} {\small Throughout this paper, $\| \cdot\|_q$ stands for the usual $L^q({\bf R}^{n})$-norm. For simplicity of notation, in particular, we use $\| \cdot\|$ instead of $\| \cdot\|_2$. 
We also introduce the following weighted functional spaces.
\[L^{1,\gamma}({\bf R}^{n}) := \left\{f \in L^{1}({\bf R}^{n}) \; \bigm| \; \Vert f\Vert_{1,\gamma} := \int_{{\bf R}^{n}}(1+\vert x\vert^{\gamma})\vert f(x)\vert dx < +\infty\right\}.\]
One denotes the Fourier transform ${\cal F}_{x\to\xi}(f)(\xi)$ of $f(x)$ by 
\[{\cal F}_{x\to\xi}(f)(\xi) = \hat{f}(\xi) := \displaystyle{\int_{{\bf R}^{n}}}e^{-ix\cdot\xi}f(x)dx, \quad \xi \in {\bf R}^n,\]
as usual with $i := \sqrt{-1}$, and ${\cal F}_{\xi\to x}^{-1}$ expresses its inverse Fourier transform. Finally, we denote the surface area of the $n$-dimensional unit ball by $\omega_{n} := \displaystyle{\int_{\vert\omega\vert = 1}}d\omega$, and we set $x = (x_{1}, x_{2}, \cdots, x_{n}) \in {\bf R}^{n}$. }

The paper is organized as follows. In Section 2 we try to get $L^{2}$-bound of solutions via known method due to \cite{IM}. In Section 3 we derive the lower bound estimates of the $L^{2}$-norm of solutions, and in Section 4 we obtain the upper bound estimates of the $L^{2}$-norm of solutions, and by combining the results obtained in Sections $3$ and $4$ one can prove Theorems 1.1 and 1.2. Finally, we prove Theorems 1.3, 1.4 and 1.5 in Section 5.\\

\section{$L^{2}$-upper bound estimates: higher dimensional case}

In this section, we introduce a device to derive $L^{2}$-estimates of solutions to problem (1.1)-(1.2) by relying on the modified Morawetz method developed newly in \cite{IM}. For this we prepare the following Hardy type inequality (see Davies-Hinz \cite{DH}). The basic concept in \cite{IM} is that the Hardy type inequality implies the $L^{2}$-upper bound estimate of the solution in several types of linear evolution equations. 
\begin{lem}\label{lem2-1} Let $n \geq 5$. Then there exists a constant $C^{*} > 0$ such that
\[\int_{{\bf R}^{n}}\frac{\vert w(x)\vert^{2}}{\vert x\vert^{4}}dx \leq C^{*}\int_{{\bf R}^{n}}\vert\Delta w(x)\vert^{2}dx\]
for all $w \in H^{2}({\bf R}^{n})$.
\end{lem}
As in the idea \cite{IM}, for the solution $u(t,x)$ to problem (1.1)-(1.2) we set
\[v(t,x) := \int_{0}^{t}u(s,x)ds.\]
Then, the function $v(t,x)$ satisfies the following equation and initial data. 
\begin{align}
& v_{tt} +\Delta^{2}v = u_{1},\ \ \ (t,x)\in (0,\infty)\times {\bf R}^{n},\label{eqn-1}\\
& v(0,x) = 0, \quad  v_{t}(0,x)= u_{0}(x),\ x\in{\bf R}^{n}.\label{initial-1}
\end{align}
Multiplying the both sides of \eqref{eqn-1}-\eqref{initial-1} by $v_{t}$, and integrating it over $[0,t]\times {\bf R}^{n}$ one can get the following energy equality such that
\begin{equation}\label{2-1}
\frac{1}{2}\Vert v_{t}(t,\cdot)\Vert^{2} + \frac{1}{2}\Vert\Delta v(t,\cdot)\Vert^{2} = \frac{1}{2}\Vert v_{0}\Vert^{2} + \int_{{\bf R}^{n}}u_{1}(x)v(t,x)dx.
\end{equation} 
Now we are in a position to use Lemma \ref{lem2-1} to proceed the following computations based on the Schwarz inequality:
\[\left\vert\int_{{\bf R}^{n}}u_{1}(x)v(t,x)dx\right\vert \leq \int_{{\bf R}^{n}}\vert u_{1}(x)\vert\vert v(t,x)\vert dx = \int_{{\bf R}^{n}}\left(\vert x\vert^{2}\vert u_{1}(x)\vert\right)\left(\frac{\vert v(t,x)\vert}{\vert x\vert^{2}}\right)dx\]
\[\leq \left(\int_{{\bf R}^{n}}\vert x\vert^{4}\vert u_{1}(x)\vert^{2}dx\right)^{1/2}\left( \int_{{\bf R}^{n}}\frac{\vert v(t,x)\vert^{2}}{\vert x\vert^{4}}dx\right)^{1/2}\]
\[\leq C_{\varepsilon}\int_{{\bf R}^{n}}\vert x\vert^{4}\vert u_{1}(x)\vert^{2}dx + \varepsilon\int_{{\bf R}^{n}}\frac{\vert v(t,x)\vert^{2}}{\vert x\vert^{4}}dx\]
\begin{equation}\label{2-2}
\leq C_{\varepsilon}\int_{{\bf R}^{n}}\vert x\vert^{4}\vert u_{1}(x)\vert^{2}dx + \varepsilon C^{*}\int_{{\bf R}^{n}} \vert\Delta v(t,x)\vert^{2}dx
\end{equation}
with some constant $C_{\varepsilon} > 0$, which depends on each $\varepsilon > 0$. Thus, from \eqref{2-1} and \eqref{2-2} it follows that
\[\frac{1}{2}\Vert v_{t}(t,\cdot)\Vert^{2} + \frac{1}{2}\Vert\Delta v(t,\cdot)\Vert^{2} \leq \frac{1}{2}\Vert u_{0}\Vert^{2} + C_{\varepsilon}\int_{{\bf R}^{n}}\vert x\vert^{4}\vert u_{1}(x)\vert^{2}dx + \varepsilon C^{*}\int_{{\bf R}^{n}} \vert\Delta v(t,x)\vert^{2}dx,\]
which implies
\[\frac{1}{2}\Vert v_{t}(t,\cdot)\Vert^{2} + \left(\frac{1}{2}-\varepsilon C^{*}\right)\Vert\Delta v(t,\cdot)\Vert^{2} \leq \frac{1}{2}\Vert u_{0}\Vert^{2} + C_{\varepsilon}\int_{{\bf R}^{n}}\vert x\vert^{4}\vert u_{1}(x)\vert^{2}dx\]
provided that
\[\int_{{\bf R}^{n}}\vert x\vert^{4}\vert u_{1}(x)\vert^{2}dx < +\infty.\]
By choosing $\varepsilon > 0$ small enough one can arrive at the crucial $L^{2}$-estimate by means of the multiplier method only because of $v_{t}(t,x) = u(t,x)$. The result below compensates the unknown information on $n \geq 5$. 
\begin{pro}\label{pro1} Let $n \geq 5$, and $[u_{0},u_{1}] \in H^{2}({\bf R}^{n})\times L^{2}({\bf R}^{n})$. Then, there exists a constant $C > 0$ such that the solution $u(t,x)$ to problem \eqref{eqn}-\eqref{initial} satisfies  
\[\Vert u(t,\cdot)\Vert ^{2} \leq C\left(\Vert u_{0}\Vert^{2} + \int_{{\bf R}^{n}}\vert x\vert^{4}\vert u_{1}(x)\vert^{2}dx\right) \quad t \geq 0,\]
provided additionally that
\begin{equation}\label{1192}
\int_{{\bf R}^{n}}\vert x\vert^{4}\vert u_{1}(x)\vert^{2}dx < +\infty.
\end{equation}
\end{pro}
\begin{rem}{\rm An essential part of this argument above has already been developed previously in \cite{LC} for the plate equation with damping and the rotational inertia term by employing the idea due to \cite{IM}. This part is just a trial to observe the case of $n \geq 5$.}
\end{rem}
\begin{rem}\label{rm1941}{\rm Note that in the case of $4 \leq n$, the condition \eqref{1192} together with $u_{1} \in L^{2}({\bf R}^{n})$ does not necessarily imply $u_{1} \in L^{1}({\bf R}^{n})$.}

\end{rem}

So, one can apply the method due to \cite{IM} to get the $L^{2}$-bound in the case when $n \geq 5$, however, this method decisively depends on the existence of the Hardy type inequality. So, the results for $n = 1,2,3,4$ of Theorems 1.1 and 1.2 are quite important since those are independent of any such inequalities. We have to derive $L^{2}$-bound estimates from the equation itself.\\
One can proceed the similar argument by using the following inequality in place of Lemma \ref{lem2-1} (see \cite{LC} for its detail).
\begin{lem}\label{lem2-2} Let $n \geq 5$. Then there exists a constant $C^{*} > 0$ such that
\[\Vert w\Vert_{L^{\frac{2n}{n-4}}} \leq C^{*}\Vert\Delta w\Vert\]
for all $w \in H^{2}({\bf R}^{n})$.
\end{lem}


\section{$L^{2}$-lower bound estimates}

In this section, we derive the lower bound estimates of the $t$-function $\Vert u(t,\cdot)\Vert $ based on the Plancherel Theorem combined with low frequency estimates (cf. \cite{I-14}). For the proof, it suffices to assume $[u_{0},u_{1}] \in C_{0}^{\infty}({\bf R}^{n}) \times C_{0}^{\infty}({\bf R}^{n})$ to have sufficient regularity of the solution $u(t,x)$ because the density argument can be applied in the final estimates. \\ 
We first set
\begin{equation}\label{i-2}
L := \sup_{\theta \ne 0}\left\vert \frac{\sin\theta}{\theta}\right\vert < +\infty.
\end{equation}
On the other hand, because of 
\[\lim_{\theta \to +0}\frac{\sin\theta}{\theta} = 1,\]
there exists a real number $\delta_{0} \in (0,1)$ such that
\begin{equation}\label{i-3}
\left\vert \frac{\sin\theta}{\theta}\right\vert \geq \frac{1}{2}
\end{equation}
for all $\theta \in (0,\delta_{0}^{2}]$. The following fundamental inequality will be useful, too:
\begin{equation}\label{i-4}
\vert a + b\vert^{2} \geq \frac{1}{2}\vert a\vert^{2} - \vert b\vert^{2} 
\end{equation}
for all $a, b \in {\bf C}$.\\
In order to get the lower bound estimate for the quantity $\Vert u(t,\cdot)\Vert$, it suffices to treat $\Vert w(t,\cdot)\Vert$ with $w(t,\xi) := {\cal F}_{x \to \xi}(u(t,\cdot))(\xi) = \hat{u}(t,\xi)$ because of the Plancherel Theorem.\\

Now we decompose the quantity $\Vert w(t,\cdot)\Vert$ as follows: for each $n \geq 1$
\begin{equation}\label{i-5}
\Vert w(t,\cdot)\Vert^{2} = \left(\int_{\vert\xi\vert \leq \frac{\delta_{0}}{\sqrt{t}}} + \int_{\vert\xi\vert \geq \frac{\delta_{0}}{\sqrt{t}}}\right)\vert w(t,\xi)\vert^{2}d\xi = I_{low}^{(n)}(t) + I_{high}^{(n)}(t).
\end{equation}
Here we have just chosen $t > 0$ large enough such that
\[\frac{\delta_{0}}{\sqrt{t}} \leq 1.\]

By the way, in the Fourier space ${\bf R}_{\xi}^{n}$ the problem \eqref{eqn}-\eqref{initial} and its solution $u(t,x)$ can be transformed into the following ODE with parameter $\xi \in {\bf R}_{\xi}^{n}$:
\begin{align}
& w_{tt} + |\xi|^4 w = 0,\ \ \ t>0,\quad \xi \in {\bf R}_{\xi}^{n},\label{eqnfourier}\\
& w(0,\xi)= w_0(\xi), \quad  w_{t}(0,\xi)= w_{1}(\xi),\ \ \ \xi \in{\bf R}^{n} ,\label{initialfourier}
\end{align}
where $w_{0}(\xi) := \hat{u}_0(\xi)$ and $w_{1}(\xi) := \hat{u}_1(\xi)$. In addition, one can solve the problem \eqref{eqnfourier}-\eqref{initialfourier} as follows:
\begin{equation}\label{i-6}
w(t,\xi) = \frac{\sin(t\vert\xi\vert^{2})}{\vert\xi\vert^{2}}w_{1}(\xi) + \cos(t\vert\xi\vert^{2})w_{0}(\xi).
\end{equation}

Now, let us derive the lower bound estimates for $I_{low}^{(n)}(t)$ because one has $\Vert w(t,\cdot)\Vert^{2} \geq I_{low}^{(n)}(t)$. One relies on a device coming from an idea in \cite{I-14}.  Indeed, it follows from \eqref{i-6} and \eqref{i-4} that
\[I_{low}^{(n)}(t) = \int_{\vert\xi\vert \leq \frac{\delta_{0}}{\sqrt{t}}}\left\vert\frac{\sin(t\vert\xi\vert^{2})}{\vert\xi\vert^{2}}w_{1}(\xi) + \cos(t\vert\xi\vert^{2})w_{0}(\xi)\right\vert^{2} d\xi\]
\[\geq \frac{1}{2}\int_{\vert\xi\vert \leq \frac{\delta_{0}}{\sqrt{t}}}\frac{\sin^{2}(t\vert\xi\vert^{2})}{\vert\xi\vert^{4}}\vert w_{1}(\xi)\vert^{2}d\xi -  \int_{\vert\xi\vert \leq \frac{\delta_{0}}{\sqrt{t}}}\cos^{2}(t\vert\xi\vert^{2})\vert w_{0}(\xi)\vert^{2}d\xi\]
\begin{equation}\label{i-7}
=: \frac{1}{2}J_{1}(t) - J_{2}(t).
\end{equation}
Let us first estimate $J_{1}(t)$ by using the decomposition of the initial data $w_{1}(\xi)$ in the Fourier space: 
\[w_{1}(\xi) = P + (A(\xi)-iB(\xi)),\quad \xi \in {\bf R}_{\xi}^{n},\]
where
\[P := \int_{{\bf R}^{n}}u_{1}(x)dx,\]
\[A(\xi) := \int_{{\bf R}^{n}}(\cos(x\xi)-1)u_{1}(x)dx, \quad B(\xi) := \int_{{\bf R}^{n}}\sin(x\xi)u_{1}(x)dx.\]
It is known (see \cite{I-04}) that there is a constant $M > 0$ such that
\begin{equation}\label{i-8}
\vert A(\xi)-iB(\xi)\vert \leq M\vert\xi\vert\Vert u_{1}\Vert_{L^{1,1}},\quad \xi \in {\bf R}_{\xi}^{n},
\end{equation}
in case of $u_{1} \in L^{1,1}({\bf R}^{n})$. Then, from \eqref{i-4} we see that
\[J_{1}(t) = \int_{\vert\xi\vert \leq \frac{\delta_{0}}{\sqrt{t}}}\frac{\sin^{2}(t\vert\xi\vert^{2})}{\vert\xi\vert^{4}}\vert w_{1}(\xi)\vert^{2}d\xi\]
\[\geq \frac{P^{2}}{2}\int_{\vert\xi\vert \leq \frac{\delta_{0}}{\sqrt{t}}}\frac{\sin^{2}(t\vert\xi\vert^{2})}{\vert\xi\vert^{4}}d\xi - \int_{\vert\xi\vert \leq \frac{\delta_{0}}{\sqrt{t}}}\vert A(\xi)-iB(\xi)\vert^{2}\frac{\sin^{2}(t\vert\xi\vert^{2})}{\vert\xi\vert^{4}}d\xi\]
\begin{equation}\label{i-9}
= \frac{P^{2}}{2}K_{1}(t)-K_{2}(t).
\end{equation}
$K_{2}(t)$ can be estimated from above by using and \eqref{i-2} and \eqref{i-8}:
\[K_{2}(t) \leq M^{2}\Vert u_{1}\Vert_{L^{1,1}}^{2}\omega_{n}\int_{0}^{\frac{\delta_{0}}{\sqrt{t}}} t^{2}\left(\frac{\sin(t\vert\xi\vert^{2})}{t\vert\xi\vert^{2}}\right)^{2}\vert\xi\vert^{2}d\xi\]
\[\leq M^{2}\Vert u_{1}\Vert_{L^{1,1}}^{2}L^{2}t^{2}\int_{\vert\xi\vert \leq \frac{\delta_{0}}{\sqrt{t}}}\vert\xi\vert^{2}d\xi = M^{2}\Vert u_{1}\Vert_{L^{1,1}}^{2}\omega_{n}L^{2}t^{2}\int_{0}^{\frac{\delta_{0}}{\sqrt{t}}}r^{n+1}dr\]
\begin{equation}\label{i-10}
= \frac{L^{2}}{n+2}M^{2}\omega_{n}\delta_{0}^{n+2}\Vert u_{1}\Vert_{L^{1,1}}^{2} t^{1-\frac{n}{2}},\quad t \gg 1.
\end{equation}
On the other hand, one can obtain the lower bound estimate for $K_{1}(t)$ because of \eqref{i-3}:
\[K_{1}(t) = t^{2}\int_{\vert\xi\vert \leq \frac{\delta_{0}}{\sqrt{t}}}\frac{\sin^{2}(t\vert\xi\vert^{2})}{(t\vert\xi\vert^{2})^{2}}d\xi \geq \frac{t^{2}}{4}\int_{\vert\xi\vert \leq \frac{\delta_{0}}{\sqrt{t}}}d\xi \]
\begin{equation}\label{i-11}
= \frac{\omega_{n}\delta_{0}^{n}}{4n}t^{2-\frac{n}{2}}, \quad t \gg 1.
\end{equation}
Therefore from \eqref{i-9}, \eqref{i-10} and \eqref{i-11} one can get the estimate from below for $J_{1}(t)$:
\begin{equation}\label{i-12}
J_{1}(t) \geq \frac{P^{2}}{2}\frac{\omega_{n}\delta_{0}^{n}}{4n}t^{2-\frac{n}{2}}-\frac{L^{2}}{n+2}M^{2}\omega_{n}\delta_{0}^{n+2}\Vert u_{1}\Vert_{L^{1,1}}^{2} t^{1-\frac{n}{2}} \quad t \gg 1.
\end{equation}
Since the upper bound estimate of $J_{2}(t)$ can be easily obtained as follows: 
\[J_{2}(t) \leq \int_{\vert\xi\vert \leq \frac{\delta_{0}}{\sqrt{t}}}\vert w_{0}(\xi)\vert^{2}d\xi \leq \Vert u_{0}\Vert_{L^{1}}^{2}\omega_{n}\int_{0}^{\frac{\delta_{0}}{\sqrt{t}}}r^{n-1}dr = \frac{\omega_{n}\delta_{0}^{n}}{n}\Vert u_{0}\Vert_{L^{1}}^{2}t^{-\frac{n}{2}},\]
because of \eqref{i-12}, \eqref{i-5} and \eqref{i-7} one has just arrived at the following lower bound estimate for $\Vert w(t,\cdot)\Vert$:
\begin{equation}\label{i-13}
\Vert w(t,\cdot)\Vert^{2} \geq I_{low}^{(n)}(t) \geq  \frac{P^{2}}{4}\frac{\omega_{n}\delta_{0}^{n}}{4n}t^{2-\frac{n}{2}}-\frac{M^{2}}{n}\omega_{n}\delta_{0}^{n}\Vert u_{1}\Vert_{L^{1,1}}^{2} t^{1-\frac{n}{2}} -\frac{\omega_{n}\delta_{0}^{n}}{n}\Vert u_{0}\Vert_{L^{1}}^{2}t^{-\frac{n}{2}}\quad t \gg 1.
\end{equation}
Therefore,  there is a positive real number $t_{0}$ such that 
\begin{equation}\label{i-13-1}
\Vert w(t,\cdot)\Vert^{2} \geq \frac{P^{2}}{32n}\omega_{n}\delta_{0}^{n} t^{2-\frac{n}{2}}
\end{equation}
for all $t \geq t_{0}$ and all $n \in {\bf N}$. It should be mentioned that $t_{0} > 0$ depends on $n$ and the quantities $\Vert u_{1}\Vert_{L^{1.1}}$ and $\Vert u_{0}\Vert_{L^{1}}$. From \eqref{i-13-1} with $n \geq 1$ one has the following lemma. 
\begin{lem}\label{lem1}Let $n \geq 1$, and $[u_{0},u_{1}] \in L^{1}({\bf R}^{n}) \times L^{1,1}({\bf R}^{n})$. Then, it holds that
\[\Vert w(t,\cdot)\Vert^{2} \geq CP^{2}t^{2-\frac{n}{2}}, \quad t \gg 1.\]
\end{lem}
Note that the results for $n = 4$ in Lemma \ref{lem1} seems to be weak because it shows just the non-decay property. We expect $\log t$-blowup as in wave equation case for $n = 2$ (\cite{Ikehata}). For this purpose we improve the statement for $n = 4$ by relying on a similar argument as in \cite{Ikehata} with a trick function $e^{-r^{2}}$. Indeed, it follows from \eqref{i-4} and a similar argument to the previous one that 
\[\Vert w(t,\cdot)\Vert^{2} \geq \frac{1}{2}\int_{{\bf R}^{4}}\frac{\sin^{2}(tr^{2})}{r^{4}}\vert w_{1}(\xi)\vert^{2}d\xi -\int_{{\bf R}^{4}}\cos^{2}(tr^{2})\vert w_{0}(\xi)\vert^{2}d\xi\]
\[\geq \frac{1}{2}\int_{{\bf R}^{4}}e^{-r^{2}}\frac{\sin^{2}(tr^{2})}{r^{4}}e^{r^{2}}\vert P+(A(\xi)-iB(\xi))\vert^{2}d\xi- \Vert u_{0}\Vert^{2}\]
\[\geq \frac{1}{4}P^{2}\int_{{\bf R}^{4}}e^{-r^{2}}\frac{\sin^{2}(tr^{2})}{r^{4}}d\xi - \frac{M^{2}\Vert u_{1}\Vert_{L^{1,1}}^{2}}{2}\int_{{\bf R}^{4}}e^{-r^{2}}\frac{\sin^{2}(tr^{2})}{r^{4}}r^{2}d\xi- \Vert u_{0}\Vert^{2}\]
\[\geq \frac{1}{4}P^{2}\int_{{\bf R}^{4}}e^{-r^{2}}\frac{\sin^{2}(tr^{2})}{r^{4}}d\xi - \frac{1}{2}M^{2}\Vert u_{1}\Vert_{L^{1,1}}^{2}\int_{{\bf R}^{4}}e^{-r^{2}}\frac{1}{r^{2}}d\xi -  \Vert u_{0}\Vert^{2}\]
\begin{equation}\label{ike-51}
=:  \frac{1}{4}P^{2}U(t) -\frac{1}{4}M^{2}\Vert u_{1}\Vert_{L^{1,1}}^{2}\omega_{4} - \Vert u_{0}\Vert^{2}.
\end{equation}
Now, we apply an useful idea coming from \cite{IO}. For this purpose, we set
\[\theta_{j} := \sqrt{(\frac{1}{4}+j)\frac{\pi}{t}}, \quad \tau_{j} := \sqrt{(\frac{3}{4}+j)\frac{\pi}{t}}\quad (j = 0,1,2,\cdots),\]
and choose $t > 1$ large enough such that $\theta_{0} = \sqrt{\frac{\pi}{4t}} < 1$. Then, because of the fact that $$\vert\sin(tr^{2})\vert \geq \frac{1}{\sqrt{2}}$$
for $r \in [\theta_{j},\tau_{j}]$ and $j = 0,1,2,\cdots$, one can proceed to estimate:
\begin{equation}\label{ike-53}
U(t) = \int_{{\bf R}^{4}}e^{-r^{2}}\frac{\sin^{2}(tr^{2})}{r^{4}}d\xi \geq \frac{1}{2}\sum_{j= 0}^{\infty}\int_{\theta_{j} \leq r \leq \tau_{j}}\frac{e^{-r^{2}}}{r^{4}}d\xi = \frac{\omega_{4}}{2}\left(\sum_{j= 0}^{\infty}\int_{\theta_{j}}^{\tau_{j}}\frac{e^{-r^{2}}}{r}dr\right)
\end{equation}
\begin{equation}\label{ike-54}
\geq \frac{\omega_{4}}{2}\left(\frac{1}{2}\int_{\theta_{0}}^{\infty}\frac{e^{-r^{2}}}{r}dr\right) \geq \frac{\omega_{4}}{4}\int_{\theta_{0}}^{1}\frac{e^{-r^{2}}}{r}dr
\end{equation}
\begin{equation}\label{ike-52}
\geq \omega_{4}\frac{e^{-1}}{4}\int_{\theta_{0}}^{1}\frac{1}{r}dr = \omega_{4}\frac{e^{-1}}{8}(\log t +\log 4 -\log \pi),
\end{equation}
where in the inequality from \eqref{ike-53} to \eqref{ike-54} one has just used the monotone decreasing property of the function $r \mapsto \frac{e^{-r^{2}}}{r}$, and the fact that the length
$$\tau_{j}-\theta_{j} = \frac{\pi}{2t}\frac{1}{\sqrt{(\frac{3}{4}+j)\frac{\pi}{t}}+\sqrt{(\frac{1}{4}+j)\frac{\pi}{t}} }$$
is decreasing to $0$ as $j \to \infty$ for each fixed $t > 1$. Therefore, by \eqref{ike-51} and \eqref{ike-52} one has the following estimates for $n = 4$ for large $t > 1$.
\begin{lem}\label{lem2}Let $n = 4$, and $[u_{0},u_{1}] \in L^{1}({\bf R}^{n}) \times L^{1,1}({\bf R}^{n})$. Then, it holds that
\[\Vert w(t,\cdot)\Vert^{2} \geq CP^{2}\log t, \quad t \gg 1.\]
\end{lem}

Proofs of the lower bound estimates of Theorems 1.1 and 1.2 are direct consequence of Lemmas \ref{lem1} and \ref{lem2}.


\section{$L^{2}$-upper bound estimates of the solution}

In this section, let us derive upper bound estimates of $\Vert u(t,\cdot)\Vert$ as $t \to \infty$ by treating the function $w(t,\xi)$ in both high and low frequency region.\\ 
As in Section 3, one again assumes $[u_{0},u_{1}] \in C_{0}^{\infty}({\bf R}^{n}) \times C_{0}^{\infty}({\bf R}^{n})$ to proceed the proof.\\

From \eqref{i-5} one first derives the upper bound estimate for $I_{low}^{(n)}(t)$ for all $n \geq 1$. Indeed, by \eqref{i-7} one has
\[I_{low}^{(n)}(t) \leq \int_{\vert\xi\vert \leq \frac{\delta_{0}}{\sqrt{t}}}\vert\frac{\sin(t\vert\xi\vert^{2})}{\vert\xi\vert^{2}}w_{1}(\xi) + \cos(t\vert\xi\vert^{2})w_{0}(\xi) \vert^{2}d\xi\]
\[\leq 2\int_{\vert\xi\vert \leq \frac{\delta_{0}}{\sqrt{t}}}\frac{\sin^{2}(t\vert\xi\vert^{2})}{\vert\xi\vert^{4}}\vert w_{1}(\xi)\vert^{2}d\xi  + 2\int_{\vert\xi\vert \leq \frac{\delta_{0}}{\sqrt{t}}}\cos^{2}(t\vert\xi\vert^{2})\vert w_{0}(\xi)\vert^{2}d\xi\]
\begin{equation}\label{i-14}
= 2L_{1}(t) + 2L_{2}(t). 
\end{equation}
It is easy to obtain the estimate for $L_{2}(t)$ as in the estimate for $J_{2}(t)$ in Section 2: 
\[L_{2}(t) \leq \int_{\vert\xi\vert \leq \frac{\delta_{0}}{\sqrt{t}}}\vert w_{0}(\xi)\vert^{2}d\xi \leq \Vert u_{0}\Vert_{L^{1}}^{2}\omega_{n}\int_{0}^{\frac{\delta_{0}}{\sqrt{t}}}r^{n-1}dr\]
\begin{equation}\label{i-15}
 = \frac{\omega_{n}\delta_{0}^{n}}{n}\Vert u_{0}\Vert_{L^{1}}^{2}t^{-\frac{n}{2}},
\end{equation}
where one has just used the fat that
\begin{equation}\label{i-17}
\int_{\vert\xi\vert \leq  \frac{\delta_{0}}{\sqrt{t}}}d\xi = \frac{\omega_{n}\delta_{0}^{n}}{n}t^{-\frac{n}{2}},\quad t \gg 1.
\end{equation}
Let us estimate $L_{1}(t)$. Indeed,  from \eqref{i-2} and \eqref{i-17} one has
\[L_{1}(t) = \int_{\vert\xi\vert \leq \frac{\delta_{0}}{\sqrt{t}}}\frac{\sin^{2}(t\vert\xi\vert^{2})}{\vert\xi\vert^{4}}\vert w_{1}(\xi)\vert^{2}d\xi \leq \Vert u_{1}\Vert_{L^{1}}^{2}t^{2}\int_{\vert\xi\vert \leq \frac{\delta_{0}}{\sqrt{t}}}\frac{\sin^{2}(t\vert\xi\vert^{2})}{(t\vert\xi\vert^{2})^{2}}d\xi\]
\begin{equation}\label{i-16}
\leq \Vert u_{1}\Vert_{L^{1}}^{2}L^{2}t^{2} \frac{\omega_{n}\delta_{0}^{n}}{n}t^{-\frac{n}{2}} = \Vert u_{1}\Vert_{L^{1}}^{2}L^{2}\frac{\omega_{n}\delta_{0}^{n}}{n}t^{2-\frac{n}{2}},\quad t \gg 1.
\end{equation}
\noindent
Thus, from \eqref{i-14}, \eqref{i-15} and \eqref{i-16} one can obtain the low-frequency estimate
\begin{equation}\label{i-low}
I_{low}^{(n)}(t) \leq C\left(\Vert u_{0}\Vert_{L^{1}}^{2}t^{-\frac{n}{2}} + \Vert u_{1}\Vert_{L^{1}}^{2}t^{2-\frac{n}{2}}\right), \quad t \gg 1,
\end{equation}
with some constant $C > 0$. 

Next, let us treat $I_{high}^{(n)}(t)$ to get the upper bound estimate. Indeed, similar to the above estimate one stars with the following inequalities.
\[I_{high}^{(n)}(t) \leq \int_{\vert\xi\vert \geq \frac{\delta_{0}}{\sqrt{t}}}\vert\frac{\sin(t\vert\xi\vert^{2})}{\vert\xi\vert^{2}}w_{1}(\xi) + \cos(t\vert\xi\vert^{2})w_{0}(\xi) \vert^{2}d\xi\]
\[\leq 2\int_{\vert\xi\vert \geq \frac{\delta_{0}}{\sqrt{t}}}\frac{\sin^{2}(t\vert\xi\vert^{2})}{\vert\xi\vert^{4}}\vert w_{1}(\xi)\vert^{2}d\xi  + 2\int_{\vert\xi\vert \geq \frac{\delta_{0}}{\sqrt{t}}}\cos^{2}(t\vert\xi\vert^{2})\vert w_{0}(\xi)\vert^{2}d\xi\]
\begin{equation}\label{i-20}
= 2N_{1}^{(n)}(t) + 2N_{2}^{(n)}(t). 
\end{equation}
It is easy to treat $N_{2}^{(n)}(t)$ as follows. This can be derived for all $n \geq 1$: 
\begin{equation}\label{i-21}
N_{2}^{(n)}(t) \leq \int_{\vert\xi\vert \geq \frac{\delta_{0}}{\sqrt{t}}}\vert w_{0}(\xi)\vert^{2}d\xi \leq \Vert u_{0}\Vert^{2}.
\end{equation}

Let us estimate $N_{1}^{(n)}(t)$. First, for all $n \geq 1$ one has 
\[N_{1}^{(n)}(t) = \int_{\vert\xi\vert \geq \frac{\delta_{0}}{t^{1/8}}}\frac{\sin^{2}(t\vert\xi\vert^{2})}{\vert\xi\vert^{4}}\vert w_{1}(\xi)\vert^{2}d\xi + \int_{\frac{\delta_{0}}{\sqrt{t}} \leq \vert\xi\vert \leq \frac{\delta_{0}}{t^{1/8}}}\frac{\sin^{2}(t\vert\xi\vert^{2})}{\vert\xi\vert^{4}}\vert w_{1}(\xi)\vert^{2}d\xi \]
\begin{equation}\label{i-22}
=: R_{1}^{(n)}(t) + R_{2}^{(n)}(t).
\end{equation}
To begin with, we treat $R_{1}^{(n)}(t)$ to get the estimate
\[R_{1}^{(n)}(t) \leq \frac{t^{1/2}}{\delta_{0}^{4}}\int_{\vert\xi\vert \geq \frac{\delta_{0}}{t^{1/8}}}\sin^{2}(t\vert\xi\vert)\vert w_{1}(\xi)\vert^{2}d\xi\]
\begin{equation}\label{i-23}
\leq \frac{t^{1/2}}{\delta_{0}^{4}}\int_{\vert\xi\vert \geq \frac{\delta_{0}}{t^{1/8}}}\vert w_{1}(\xi)\vert^{2}d\xi \leq  \frac{\sqrt{t}}{\delta_{0}^{4}}\Vert u_{1}\Vert^{2}, \quad t \gg 1.
\end{equation}
\noindent
On the other hand, in the case of $n \ne 4$, $R_{2}^{(n)}(t)$ can be estimated as follows: 
\[R_{2}^{(n)}(t) \leq \Vert u_{1}\Vert_{1}^{2}\int_{\frac{\delta_{0}}{\sqrt{t}} \leq \vert\xi\vert \leq \frac{\delta_{0}}{t^{1/8}}}\frac{1}{\vert\xi\vert^{4}}d\xi\]
\begin{equation}\label{i-25}
= \Vert u_{1}\Vert_{1}^{2}\omega_{n}\int_{\frac{\delta_{0}}{\sqrt{t}}}^{\frac{\delta_{0}}{t^{1/8}}}r^{n-5}dr = \frac{\omega_{n}}{(4-n)\delta_{0}^{4-n}}\Vert u_{1}\Vert_{1}^{2}(t^{2-\frac{n}{2}}-t^{\frac{1}{2}-\frac{n}{8}}), \quad t \gg 1.
\end{equation}
Thus, from \eqref{i-20}, \eqref{i-21}, \eqref{i-22}, \eqref{i-23} and \eqref{i-25} one has the infinite time blowup estimate for $I_{high}^{(n)}(t)$ in the case of $n = 1,2,3$:
\[
I_{high}^{(n)}(t) = C\left(\Vert u_{0}\Vert^{2} + \Vert u_{1}\Vert^{2}t^{\frac{1}{2}} + \Vert u_{1}\Vert_{L^{1}}^{2}\frac{1}{4-n}(t^{2-\frac{n}{2}}-t^{\frac{1}{2}-\frac{n}{8}})\right), \quad t \gg 1.
\]
\begin{equation}\label{i-24}
\leq C\left(\Vert u_{0}\Vert^{2} + \Vert u_{1}\Vert^{2} + \Vert u_{1}\Vert_{L^{1}}^{2}\frac{1}{4-n}\right)t^{2-\frac{n}{2}}, \quad t \gg 1.
\end{equation}

Next, let us give sharp estimates for $N_{1}^{(4)}(t)$ in the case of $n = 4$ by using a more delicate computation. For this purpose, by choosing $t > 1$ sufficiently large to get the relation $1 \leq \log t \leq t \leq t^{2}$ one has a decomposition of the integrand:
\[N_{1}^{(4)}(t) = \int_{\vert\xi\vert \geq \frac{\delta_{0}}{(\log t)^{1/4}}}\frac{\sin^{2}(t\vert\xi\vert^{2})}{\vert\xi\vert^{4}}\vert w_{1}(\xi)\vert^{2}d\xi\]
\[+ \int_{\frac{\delta_{0}}{(\log t)^{1/4}} \geq \vert\xi\vert \geq \frac{\delta_{0}}{t^{1/4}}}\frac{\sin^{2}(t\vert\xi\vert^{2})}{\vert\xi\vert^{4}}\vert w_{1}(\xi)\vert^{2}d\xi + \int_{\frac{\delta_{0}}{t^{1/4}} \geq \vert\xi\vert \geq \frac{\delta_{0}}{\sqrt{t}}}\frac{\sin^{2}(t\vert\xi\vert^{2})}{\vert\xi\vert^{4}}\vert w_{1}(\xi)\vert^{2}d\xi \]
\begin{equation}\label{i-22-10}
=: S_{1}(t) + S_{2}(t) + S_{3}(t).
\end{equation}
Let us estimate them in order.\\
First, one has
\[S_{1}(t) \leq \frac{\log t}{\delta_{0}^{4}}\int_{\vert\xi\vert \geq \frac{\delta_{0}}{(\log t)^{1/4}}}\vert w_{1}(\xi)\vert^{2}d\xi\]
\begin{equation}\label{i-22-1}
\leq  \frac{\log t}{\delta_{0}^{4}}\Vert u_{1}\Vert^{2},\quad t \gg 1.
\end{equation}
For  $S_{2}(t)$ we see that
\[ S_{2}(t) \leq \Vert u_{1}\Vert_{L^{1}}^{2}\omega_{4}\int_{\frac{\delta_{0}}{t^{1/4}}}^{\frac{\delta_{0}}{(\log t)^{1/4}}}\frac{1}{r}dr\]
\begin{equation}\label{i-22-2}
\leq  \omega_{4}\Vert u_{1}\Vert_{L^{1}}^{2}4^{-1}\left(\log t - \log(\log t)\right),\quad t \gg 1.
\end{equation}
Finally, let us treat $S_{3}(t)$ to obtain the following estimate:
\[S_{3}(t) \leq \Vert u_{1}\Vert_{L^{1}}^{2}\int_{\frac{\delta_{0}}{\sqrt{t}} \leq \vert\xi\vert \leq \frac{\delta_{0}}{t^{1/4}}}\frac{1}{r^{4}}d\xi\] 
\begin{equation}\label{i-22-3}
\leq  \omega_{4}\Vert u_{1}\Vert_{L^{1}}^{2}4^{-1}\log t,\quad t \gg 1.
\end{equation}
Thus, it follows from \eqref{i-22-10}, \eqref{i-22-1}, \eqref{i-22-2} and \eqref{i-22-3} that for large $t \gg 1$
\begin{equation}\label{i-22-4}
N_{1}^{(4)}(t) \leq C\left(\Vert u_{1}\Vert_{L^{1}}^{2} \log t + \Vert u_{1}\Vert^{2}\log t + \Vert u_{1}\Vert_{L^{1}}^{2}(\log t-\log(\log t))\right), \quad t \gg 1
\end{equation}
with some constant $C > 0$.
\noindent
Therefore, by combining \eqref{i-20}, \eqref{i-21} and \eqref{i-22-4} one can get the high-frequency estimate for $n = 4$:
\begin{equation}\label{i-high-2}
I_{high}^{(4)}(t) \leq C\left(\Vert u_{0}\Vert^{2} + \Vert u_{1}\Vert^{2} + \Vert u_{1}\Vert_{L^{1}}^{2}\right)\log t, \quad t \gg 1.
\end{equation}

Finally, it follows from \eqref{i-5}, \eqref{i-low}, \eqref{i-24} and \eqref{i-high-2} one can get the crucial upper bound estimates for each $n = 1, 2, 3, 4$.
\begin{lem}\label{lem3}Let $n = 1, 2, 3$, and $[u_{0},u_{1}] \in (L^{1}({\bf R}^{n})\cap L^{2}({\bf R}^{n})) \times (L^{1}({\bf R}^{n})\cap L^{2}({\bf R}^{n}))$. Then, it holds that
\[\Vert w(t,\cdot)\Vert^{2} \leq C(\Vert u_{0}\Vert^{2} + \Vert u_{0}\Vert_{L^{1}}^{2} + \Vert u_{1}\Vert^{2} + \Vert u_{1}\Vert_{L^{1}}^{2})t^{2-\frac{n}{2}}, \quad t \gg 1.\]
\end{lem}
\begin{lem}\label{lem4}Let $n = 4$, and $[u_{0},u_{1}] \in (L^{1}({\bf R}^{n})\cap L^{2}({\bf R}^{n})) \times (L^{1}({\bf R}^{n})\cap L^{2}({\bf R}^{n}))$. Then, it holds that
\[\Vert w(t,\cdot)\Vert^{2} \leq C(\Vert u_{0}\Vert^{2} +\Vert u_{0}\Vert_{L^{1}}^{2} + \Vert u_{1}\Vert^{2} + \Vert u_{1}\Vert_{L^{1}}^{2})\log t, \quad t \gg 1.\]
\end{lem}
\par
\vspace{0.5cm}
Finally, proofs of the upper bound estimates part of Theorems \ref{theorem1} and \ref{theorem2} are direct consequence of Lemmas \ref{lem3} and \ref{lem4}, and the Plancherel Theorem.\\
In connection with these estimates above, one can show the $L^{2}$-upper bound estimate of the solution itself independently from the Hardy type inequality (see Proposition \ref{pro1} and Remark \ref{rm1941}).
\begin{pro}\label{pro2} Let $n \geq 5$, and $[u_{0},u_{1}] \in H^{2}({\bf R}^{n}) \times (L^{2}({\bf R}^{n})\cap L^{1}({\bf R}^{n}))$. Then, there exists a constant $C > 0$ such that the solution $u(t,x)$ to problem \eqref{eqn}-\eqref{initial} satisfies  
\[\Vert u(t,\cdot)\Vert \leq C\left(\Vert u_{0}\Vert + \Vert u_{1}\Vert + \Vert u_{1}\Vert_{1} \right) \quad (t \geq 0).\]
\end{pro}
\begin{rem}{\rm A similar $L^{2}$-upper bound estimates can be derived to the solution itself of the (free) wave equation in the Euclidean space ${\bf R}^{n}$ in the case of $n \geq 3$, however, the method developed in Proposition 2.1 is widely applicable to the exterior problem and variavle coefficient cases.}
\end{rem}
{\it Proof of Proposition \ref{pro2}.} It suffices to assume the initial data $u_{j}$ ($j = 0,1$) belongs to $C_{0}^{\infty}({\bf R}^{n})$. Then, one can proceed to estimate as follows:
\[C_{n}\Vert u(t,\cdot)\Vert^{2} = \Vert w(t,\cdot)\Vert^{2} = \int_{{\bf R}^{n}}\left\vert\frac{\sin(t\vert\xi\vert^{2})}{\vert\xi\vert^{2}}w_{1}(\xi) + \cos(t\vert\xi\vert^{2})w_{0}(\xi) \right\vert^{2}d\xi\]
\[\leq 2\int_{{\bf R}^{n}}\frac{\sin^{2}(t\vert\xi\vert^{2})}{\vert\xi\vert^{4}}\vert w_{1}(\xi)\vert^{2}d\xi  + 2\int_{{\bf R}^{n}}\cos^{2}(t\vert\xi\vert^{2})\vert w_{0}(\xi)\vert^{2}d\xi\]
\[\leq  2\int_{\vert\xi\vert\leq 1}\frac{\sin^{2}(t\vert\xi\vert^{2})}{\vert\xi\vert^{4}}\vert w_{1}(\xi)\vert^{2}d\xi  +  2\int_{\vert\xi\vert\geq1}\frac{\sin^{2}(t\vert\xi\vert^{2})}{\vert\xi\vert^{4}}\vert w_{1}(\xi)\vert^{2}d\xi  + \Vert u_{0}\Vert^{2}\]
\begin{equation}\label{i-141}
= 2T_{1}(t) + 2T_{2}(t) + \Vert u_{0}\Vert^{2}. 
\end{equation}
Now, about $T_{1}(t)$ one can get
\[T_{1}(t) \leq \Vert w_{1}\Vert_{\infty}^{2}\int_{\vert\xi\vert\leq 1}\frac{1}{\vert\xi\vert^{4}}d\xi \leq \omega_{n}\Vert u_{1}\Vert_{1}^{2}\int_{0}^{1}r^{n-5}dr\]
\begin{equation}\label{i-142}
= \frac{\omega_{n}}{n-4}\Vert u_{1}\Vert_{1}^{2}, \quad (t \geq 0).
\end{equation}
For the estimate of $T_{2}(t)$ one can obtain
\[T_{2}(t) \leq \int_{\vert\xi\vert\geq1}\sin^{2}(t\vert\xi\vert^{2})\vert w_{1}(\xi)\vert^{2}d\xi \leq \int_{\vert\xi\vert\geq1}\vert w_{1}(\xi)\vert^{2}d\xi\]
\begin{equation}\label{i-143}
\leq  \Vert w_{1}\Vert^{2} = \Vert u_{1}\Vert^{2}, \quad (t \geq 0). 
\end{equation}
The statement of Proposition \ref{pro2} can be derived because of \eqref{i-141}-\eqref{i-143}.\\

\hfill
$\Box$\\


\section{Proof of Theorems 1.3, 1.4 and 1.5.}

In this section, let us prove Theorems 1.3, 1.4 and 1.5.  
The first part is devoted to the proof of Theorem 1.3.\\
{\it Proof of Theorem 1.3.}\, The whole idea comes from \cite{Ike-Aze}. By applying the spatial Fourier transform to the both sides of (1.1), again the problem (1.1)-(1.2) can be reduced to ODE with parameter $\xi$:
\begin{equation}
\hat{u}_{tt}(t,\xi) + \vert\xi\vert^{4}\hat{u}(t,\xi) = 0,\ \ \ \ (t,\xi) \in (0,\infty) \times {\bf R}_{\xi}^{n},
\end{equation}
\begin{equation}
\hat{u}(0,\xi) = \hat{u}_{0}(\xi),\quad \hat{u}_{t}(0,\xi) = \hat{u}_{1}(\xi),\,\,\,\, \xi \in {\bf R}_{\xi}^{n}. 
\end{equation}
For the solution $\hat{u}(t,\xi)$ to problem (5.1)-(5.2), one introduces an auxiliary function
$$\hat{w}(t,\xi) := \int^{t}_{0}\hat{u}(s,\xi)ds.$$
Then $\hat{w}(t,\xi)$ satisfies
\begin{equation}
\hat{w}_{tt}(t,\xi) + \vert\xi\vert^{4}\hat{w}(t,\xi) = \hat{u}_{1}(\xi),\ \ \ \ (t,\xi) \in (0,\infty) \times {\bf R}_{\xi}^{n},
\end{equation}
\begin{equation}
\hat{w}(0,\xi) = 0,\quad \hat{w}_{t}(0,\xi) = \hat{u}_{0}(\xi),\,\,\,\, \xi \in {\bf R}_{\xi}^{n}. 
\end{equation}
Now, we introduce one more auxiliary function $\hat{v}(t,\xi)$ defined on ${\bf R}_{\xi}^{n}\setminus\{0\}$ as follows:
\begin{equation}
\hat{v}(t,\xi) = \hat{w}(t,\xi) - \frac{\hat{u}_{1}(\xi)}{\vert\xi\vert^{4}}, \quad \xi \in {\bf R}_{\xi}^{n}\setminus\{0\}.
\end{equation}
Then, the function $\hat{v}(t,\xi)$ satisfies
\begin{equation}
\hat{v}_{tt}(t,\xi) + \vert\xi\vert^{4}\hat{v}(t,\xi) = 0,\ \ \ \  t > 0, \quad \xi \in  {\bf R}_{\xi}^{n}\setminus\{0\},
\end{equation}
\begin{equation}\label{1500}
\hat{v}(0,\xi) = -\frac{\hat{u}_{1}(\xi)}{\vert\xi\vert^{4}},\quad \hat{v}_{t}(0,\xi) = \hat{u}_{0}(\xi),\,\,\,\, \xi \in  {\bf R}_{\xi}^{n}\setminus\{0\}. 
\end{equation}
By multiplying both sided of (5.6) by $\overline{\hat{v}_{t}(t,\xi)}$, integrating it over $\{\vert \xi\vert \geq \delta\}$ with small $\delta > 0$, and taking real parts of the resulted equality one can gets
\[\frac{d}{dt}\int_{\vert\xi\vert\geq\delta}(\vert\hat{v}_{t}(t,\xi)\vert^{2} + \vert\xi\vert^{4}\vert\hat{v}(t,\xi)\vert^{2})d\xi = 0,\]
so that by integrating it over $[0,t]$ one has
\[\int_{\vert\xi\vert\geq\delta}(\vert\hat{v}_{t}(t,\xi)\vert^{2} + \vert\xi\vert^{4}\vert\hat{v}(t,\xi)\vert^{2})d\xi = \int_{\vert\xi\vert\geq\delta}(\vert\hat{u}_{0}(\xi)\vert^{2} + \frac{\vert\hat{u}_{1}(\xi)\vert^{2}}{\vert\xi\vert^{4}})d\xi,\]
where we have just used \eqref{1500}. Since $\hat{v}_{t}(t,\xi) = \hat{w}_{t}(t,\xi) = \hat{u}(t,\xi)$, one can arrive at

\[\int_{\vert\xi\vert\geq\delta}\vert\hat{u}(t,\xi)\vert^{2}d\xi \leq \int_{{\bf R}_{\xi}^{n}}\vert\hat{u}_{0}(\xi)\vert^{2}d\xi + \int_{\vert\xi\vert\geq\delta}\frac{\vert\hat{u}_{1}(\xi)\vert^{2}}{\vert\xi\vert^{4}}d\xi\]
\[\leq \Vert u_{0}\Vert^{2} + \int_{1 \geq \vert\xi\vert \geq \delta}\frac{\vert\hat{u}_{1}(\xi)\vert^{2}}{\vert\xi\vert^{4}}d\xi + \int_{1 \leq \vert\xi\vert}\frac{\vert\hat{u}_{1}(\xi)\vert^{2}}{\vert\xi\vert^{4}}d\xi\]
\begin{equation}\label{100}
\leq \Vert u_{0}\Vert^{2} + \int_{1 \geq \vert\xi\vert\geq\delta}\frac{\vert\hat{u}_{1}(\xi)\vert^{2}}{\vert\xi\vert^{4}}d\xi + \Vert u_{1}\Vert^{2}.
\end{equation}
By applying Lemma 3.1 of \cite[page 879]{I-04} to $\hat{u}_{1}(\xi)$ one has
\begin{equation}\label{101}
\vert \hat{u}_{1}(\xi)\vert  \leq C\vert\xi\vert^{\gamma}\Vert u_{1}\Vert_{1,\gamma},
\end{equation}
because of the assumption $P_{1} = \displaystyle{\int_{{\bf R}^{n}}}u_{1}(x)dx = 0$, where $\gamma \in [0,1]$. 
By substituting \eqref{101} into \eqref{100} one has arrived at the estimate
\[\int_{\vert\xi\vert\geq\delta}\vert\hat{u}(t,\xi)\vert^{2} \leq \Vert u_{0}\Vert^{2} + \Vert u_{1}\Vert^{2} + C\Vert u_{1}\Vert_{1,\gamma}^{2}\int_{1 \geq \vert\xi\vert\geq \delta}\vert\xi\vert^{2\gamma-4}d\xi\]
\begin{equation}\label{102}
\leq \Vert u_{0}\Vert^{2} + \Vert u_{1}\Vert^{2} + C\Vert u_{1}\Vert_{1,\gamma}^{2}\frac{\omega_{n}}{2\gamma + n-4},
\end{equation}
provided that $2\gamma+n-4 > 0$ for $\gamma \in (0,1]$ in the case of $n = 4$, and for $\gamma \in (\frac{1}{2},1]$ in the case of $n = 3$. Note that the constant $C > 0$ does not depend on any small $\delta > 0$. By letting $\delta \to +0$ in \eqref{102}, one has the desired estimates via the Plancherel Theorem.
$\hfill\Box$

In order to prove Theorem 1.4, we use the following lemma, which is a direct consequence of Lemma 5.1 with $\gamma \in (1,2]$ of \cite{IMichi}.
\begin{lem}\label{103}\,Let $n \geq 1$, and for $\gamma \in (1,2]$, suppose that  $f \in L^{1,\gamma}({\bf R}^{n})$ satisfies
\[\int_{{\bf R}^{n}}f(x)dx = 0, \qquad \int_{{\bf R}^{n}}x_{j}f(x)dx = 0\quad (j = 1,2,\cdots,n).\]
Then, it holds that
\[\vert \hat{f}(\xi)\vert \leq C\vert\xi\vert^{\gamma}\Vert f\Vert_{1,\gamma},\quad \xi \in {\bf R}^{n}\]
with some constant $C = C_{\gamma} > $.
\end{lem}
Based on Lemma \ref{103}, let us prove Theorem 1.4. \\

{\it Proof of Theorem 1.4.}\,The proof is just a slight modification of those for Theorem 1.3. Indeed, we use Lemma \ref{103} with $f(x) := u_{1}(x)$ and $\gamma \in (1,2]$ in place of \eqref{102} with $\gamma \in [0,1]$. Then, \eqref{102} above can be interpreted as follows:
the positivity $2\gamma+n-4 > 0$ again holds for $\gamma \in (1,2]$ in the case of $n = 2$, and for $\gamma \in (\frac{3}{2},2]$ in the case of $n = 1$. This implies the desired results.
$\hfill\Box$\\

Finally, let us prove Theorem 1.5 basing on the following equality which is a direct consequence of \cite[Lemma 5.1]{IMichi} with $f(x) := u_{1}(x)$: 
\begin{equation}\label{104}
w_{1}(\xi) = P + i\xi P_{1} + E(\xi),
\end{equation}
where
\[P := \int_{{\bf R}}u_{1}(x)dx, \quad P_{1} := \int_{{\bf R}}x u_{1}(x)dx,\]
and the error term $E(\xi)$ satisfies
\begin{equation}\label{105}
\vert E(\xi)\vert \leq C\vert\xi\vert^{2}\Vert u_{1}\Vert_{1,2}, \quad \xi \in {\bf R}.
\end{equation}

{\it Proof of (1) of Theorem 1.5.}\,For our purpose it suffices to obtain the lower bound estimate on $J_{1}(t)$ in \eqref{i-7}. Indeed, from \eqref{104} and \eqref{i-4} one has
\[J_{1}(t) = \int_{ \vert\xi\vert\leq\frac{\delta_{0}}{\sqrt{t}}}\frac{\sin^{2}(t\vert\xi\vert^{2})}{\vert\xi\vert^{4}}\vert P + i\xi P_{1} + E(\xi)\vert^{2}d\xi\] 
\[\geq \frac{1}{2}\int_{\vert\xi\vert\leq\frac{\delta_{0}}{\sqrt{t}}}\frac{\sin^{2}(t\vert\xi\vert^{2})}{\vert\xi\vert^{4}}(P^{2} + \vert\xi\vert^{2}P_{1}^{2})d\xi - \int_{\vert\xi\vert\leq\frac{\delta_{0}}{\sqrt{t}}}\frac{\sin^{2}(t\vert\xi\vert^{2})}{\vert\xi\vert^{4}}\vert E(\xi)\vert^{2}d\xi\]
\[\geq \frac{P^{2}}{2}\int_{\vert\xi\vert\leq\frac{\delta_{0}}{\sqrt{t}}}\frac{\sin^{2}(t\vert\xi\vert^{2})}{\vert\xi\vert^{4}}d\xi + \frac{P_{1}^{2}}{2}\int_{\vert\xi\vert\leq\frac{\delta_{0}}{\sqrt{t}}}\frac{\sin^{2}(t\vert\xi\vert^{2})}{\vert\xi\vert^{2}}d\xi - \int_{\vert\xi\vert\leq\frac{\delta_{0}}{\sqrt{t}}}\frac{\sin^{2}(t\vert\xi\vert^{2})}{\vert\xi\vert^{4}}\vert E(\xi)\vert^{2}d\xi\]
\begin{equation}\label{106}
=: \frac{P^{2}}{2}K_{1}(t) + \frac{P_{1}^{2}}{2}K_{3}(t) - K_{4}(t),
\end{equation}
where $K_{1}(t)$ is the same function defined in \eqref{i-9}, and 
\[K_{3}(t) := \int_{\vert\xi\vert\leq\frac{\delta_{0}}{\sqrt{t}}}\frac{\sin^{2}(t\vert\xi\vert^{2})}{\vert\xi\vert^{2}}d\xi,\]
\[K_{4}(t) :=\int_{\vert\xi\vert\leq\frac{\delta_{0}}{\sqrt{t}}}\frac{\sin^{2}(t\vert\xi\vert^{2})}{\vert\xi\vert^{4}}\vert E(\xi)\vert^{2}d\xi.\]
First, it is easy to see by using \eqref{105}
\begin{equation}\label{107}
K_{4}(t) \leq C\Vert u_{1}\Vert_{1,2}^{2}t^{-\frac{n}{2}}, \quad (n = 1)
\end{equation}
with some constant $C > 0$. While,
\[K_{3}(t) \geq \frac{t}{\delta_{0}^{2}}\omega_{1}\int_{0}^{\frac{\delta_{0}}{\sqrt{t}}}\sin^{2}(tr^{2})dr\]
\begin{equation}\label{108}
= \left(\frac{\omega_{1}}{\delta_{0}^{2}}\int_{0}^{\delta_{0}}\sin^{2}(\sigma^{2})d\sigma\right)t^{\frac{1}{2}} \quad (t \gg 1).
\end{equation}
On the lower bund of $K_{1}(t)$ one has already derived it in \eqref{i-11} wth $n = 1$. Thus, by combining it with \eqref{106}-\eqref{108} one has the desired estimates.

$\hfill\Box$\\

{\it Proof of (2) of Theorem 1.5.} Let $n = 2$. The following result, once more, is a direct consequence of \cite[Lemma 5.1]{IMichi} with $f(x) := u_{1}(x)$: 
\begin{equation}\label{104}
w_{1}(\xi) = P + i\xi \cdot {\bf P}_{1} + E(\xi),
\end{equation}
where
\[P := \int_{{\bf R^{2}}}u_{1}(x)dx, \quad {\bf P}_{1} := (p_{1},p_{2}),\]
\[p_{1} := \int_{{\bf R^{2}}}x_{1}u_{1}(x)dx,\quad p_{2} := \int_{{\bf R^{2}}}x_{2}u_{1}(x)dx,\]
and the error term $E(\xi)$ satisfies
\begin{equation}\label{105}
\vert E(\xi)\vert \leq C\vert\xi\vert^{2}\Vert u_{1}\Vert_{1,2}, \quad \xi \in {\bf R}^{2}.
\end{equation}
Similarly to the estimates developed in \eqref{ike-51}, by using \eqref{105} one can proceed the computation as follows:  
\[\Vert w(t,\cdot)\Vert^{2} \geq \frac{1}{2}\int_{{\bf R}^{2}}\frac{\sin^{2}(tr^{2})}{r^{4}}\vert w_{1}(\xi)\vert^{2}d\xi -\int_{{\bf R}^{2}}\cos^{2}(tr^{2})\vert w_{0}(\xi)\vert^{2}d\xi\]
\[\geq \frac{1}{2}\int_{{\bf R}^{2}}e^{-r^{2}}\frac{\sin^{2}(tr^{2})}{r^{4}}\vert P+ i\xi \cdot {\bf P}_{1} + E(\xi)\vert^{2}d\xi- \Vert u_{0}\Vert^{2}\]
\[\geq \frac{1}{4}P^{2}\int_{{\bf R}^{2}}e^{-r^{2}}\frac{\sin^{2}(tr^{2})}{r^{4}}d\xi + \frac{1}{4}K_{5}(t) - C^{2}\Vert u_{1}\Vert_{L^{1,2}}^{2}\int_{{\bf R}^{4}}e^{-r^{2}}\frac{\sin^{2}(tr^{2})}{r^{4}}r^{4}d\xi- \Vert u_{0}\Vert^{2}\]
\begin{equation}\label{106}
=:  \frac{1}{4}P^{2}U(t) + \frac{1}{4}K_{5}(t) -C\Vert u_{1}\Vert_{L^{1,2}}^{2}\omega_{2} - \Vert u_{0}\Vert^{2},
\end{equation}
where
\[K_{5}(t) := \int_{{\bf R}^{2}}e^{-r^{2}}\frac{\sin^{2}(tr^{2})}{r^{4}}\vert\xi\cdot{\bf P}_{1}\vert^{2}d\xi.\]

First, by similar computations as in \eqref{ike-53}, \eqref{ike-54} and \eqref{ike-52} one can get
\begin{equation}\label{107}
U(t) \geq Ct,\quad t \gg 1,
\end{equation}
with some constant $C > 0$. So, it suffices to estimate $K_{5}(t)$. In order to estimate $K_{5}(t)$, we employ an idea coming from \cite{IO}. For this, we set
\[K := \left\{ \xi \in {\bf R}^{2}\,:\,\left\vert \frac{\xi}{\vert\xi\vert}\cdot\frac{{\bf P}_{1}}{\vert {\bf P}_{1}\vert}\right\vert \geq \frac{1}{2}\right\}.\]
Then, one can estimate as follows:
\[K_{5}(t) := \int_{{\bf R}^{2}}e^{-r^{2}}\frac{\sin^{2}(tr^{2})}{r^{4}}\vert\xi\cdot{\bf P}_{1}\vert^{2}d\xi\]
\[\geq \vert {\bf P}_{1}\vert^{2}\int_{K}e^{-r^{2}}\frac{\sin^{2}(tr^{2})}{r^{2}}\vert\frac{\xi}{\vert\xi\vert}\cdot\frac{{\bf P}_{1}}{\vert{\bf P}_{1}\vert}\vert^{2}d\xi \geq \frac{\vert {\bf P}_{1}\vert^{2}}{4}\int_{K}e^{-r^{2}}\frac{\sin^{2}(tr^{2})}{r^{2}}d\xi\]
\begin{equation}\label{108}
= C\pi\frac{\vert {\bf P}_{1}\vert^{2}}{4}\int_{0}^{\infty}e^{-r^{2}}\frac{\sin^{2}(tr^{2})}{r}dr =: C\pi\frac{\vert {\bf P}_{1}\vert^{2}}{4}K_{6}(t).
\end{equation}

Once more, as in \eqref{ike-53}, \eqref{ike-54} and \eqref{ike-52} one can estimate $K_{6}(t)$:
\[K_{6}(t) \geq \sum_{j=0}^{\infty}\int_{\theta_{j}}^{\tau_{j}}e^{-r^{2}}\frac{\sin^{2}(tr^{2})}{r}dr \geq \frac{1}{2}\sum_{j=0}^{\infty}\int_{\theta_{j}}^{\tau_{j}}e^{-r^{2}}r^{-1}dr \]
\[\geq \frac{1}{4}\int_{\theta_{0}}^{\infty}e^{-r^{2}}r^{-1}dr \geq \frac{e^{-1}}{4}\int_{\theta_{0}}^{1}r^{-1}dr\]
\begin{equation}\label{109}
 \geq C\log t,\quad t \gg 1.
\end{equation}

The desired estimate can be established from \eqref{106}-\eqref{109}.

$\hfill\Box$\\

\par
\vspace{0.5cm}
\noindent{\em Acknowledgement.}
\smallskip
The author would like to thank Wenhui Chen (Shanghai Jiao Tong University) for his fruitful discussion with me. The work of the author was supported in part by Grant-in-Aid for Scientific Research (C)20K03682  of JSPS.



\begin{thebibliography}{99}




\bibitem{D} M. D'Abbicco, G. Girardi and J. Liang, $L^{1}$-$L^{1}$ estimates for the strongly damped plate equation, J. Math. Anal. Appl. {\bf 478} (2019), 476--498.

\bibitem{LC} C. R. da Luz and R. C. Chara\~o, Asymptotic properties for a semilinear plate equation in unbounded domains, J. Hyperbolic Differ. Equ. {\bf 6} (2009), 269--294.

\bibitem{I-10} R. C. Chara\~o and R. Ikehata Asymptotic profile and optimal decay of solutions of some wave equation with logarithmic damping, Z. Angew. Math. Phys. {\bf 71}, no. 148 (2020). DOI: 10.1007/s00033-020-01373-x


\bibitem{DH} E. B. Davies and A. M. Hinz, Explicit constants for Rellich inequalities in $L^{p}(\Omega)$, Math. Z. {\bf 227} (1998), 511--523.


\bibitem{I-04} R. Ikehata, New decay estimates for linear damped wave equations and its application to nonlinear problem, Math. Meth. Appl. Sci. {\bf 27} (2004), 865-889. doi: 10.1002/mma.476.




\bibitem{I-14} R. Ikehata, Asymptotic profiles for wave equations with strong damping, J. Diff. Eqns {\bf 257} (2014), 2159-2177.

\bibitem{Ikehata} R. Ikehata, $L^{2}$-blowup estimates of the wave equation and its application to local energy decay, arXiv: 2111.0203lv2 [math.AP] 13 Nov 2021.

\bibitem{Ike-Aze} R. Ikehata, Some remarks on the local energy decay for wave equations in the whole space, Azerbaijan Math. J. {\bf 9}, no. 2, (2019), 167--182.

\bibitem{IM} R. Ikehata and T. Matsuyama, $L^{2}$-behaviour of solutions to the linear heat and wave equations in exterior domains, Sci. Math. Jpn {\bf 55} (2002), 33-44.


\bibitem{IMichi} R. Ikehata and H. Michihisa, Moment conditions and lower bounds in expanding solutions of wave equations with double sdamping terms, Asymptotic Anal. {\bf 114} (2019), 19--36.

\bibitem{IO} R. Ikehata and M. Onodera, Remarks on large time behavior of the $L^{2}$-norm of solutions to strongly damped wave equations, Differ. Integral Equ. {\bf 30} (2017), 505--520.

\bibitem{SJ} S. Jing, $L^{p}$-$L^{q}$ estimates for solutions to the damped plate equation in exterior domains, Results in Math. {\bf 18} (1990), 231--253.

\bibitem{L} S. P. Levandosky, Decay estimates for fourth order wave equations, J. Diff. Eqns {\bf 143} (1998), 360--413.



\bibitem{St-2} S. P. Levandosky and W. A. Strauss, Time decay for the nonlinear beam equation, Methods Appl. Anal. {\bf 3} (2000), 479--488.

\bibitem{lin} J. E. Lin, Local time decay for a nonlinear beam equation, Methods Appl. Anal. {\bf 11} (2004), 65--68.

\bibitem{Miao} C. Miao, A note on time decay for the nonlinear beam equation, J. Math. Anal. Appl. {\bf 314} (2006), 764--773.







\bibitem{P} H. Pecher, $L^{p}$-Absch\"atzungen und klassische L\"osungen f\"ur nichtlineare Wellengleichungen. I. Math. Z. {\bf 150} (1976), 159--183.


\bibitem{peral} J. C. Peral, $L^{p}$ estimates for the wave equation, J. Func. Anal. {\bf 36} (1980), 114--145.

\bibitem{Racke} R. Racke, Lectures on Nonlinear Evolution Equations: Initial value Problems, Vieweg Verlag, Braunscheweig, 1992.


















\end{thebibliography}
\end{document}